\documentclass{amsart}
\usepackage{graphicx}
\usepackage{amsmath}%in Winedt
\usepackage{amssymb}
\usepackage{amsthm}
\usepackage{esint}
\usepackage{xcolor}

\newtheorem{theorem}{Theorem}[section]

\theoremstyle{Corollary}

\numberwithin{equation}{section}

\begin{document}

\title{Free boundary minimal annuli in $S^2_+\times S^1$}

\author{Pak Tung Ho}
%    Address of record for the research reported here
\address{Department of Mathematics, Tamkang University, Tamsui, New Taipei City 251301, Taiwan}

\email{paktungho@yahoo.com.hk}

\author{Juncheol Pyo}

\address{Department of Mathematics and Institute of Mathematical Science, Pusan National University, Busan 46241, Korea and Korea Institute for Advanced Study, Seoul 02455, Korea}

\email{jcpyo@pusan.ac.kr}

\author{Keomkyo Seo}

\address{Department of Mathematics and Research Institute of Natural Sciences, Sookmyung Women’s University, Cheongpa-ro 47-gil 100, Yongsan-ku, Seoul, 04310, Korea}

\email{kseo@sookmyung.ac.kr}

%    General info
\subjclass[2020]{Primary 53A10; Secondary 53C42}

\date{Septembert 1th, 2024.}

\begin{abstract}
Let $M$ be a compact 3-dimensional Riemannian manifold with nonnegative Ricci curvature and a nonempty boundary $\partial M$. Fraser and Li \cite{Fraser&Li} established a compactness theorem for the space of compact, properly embedded minimal surfaces of fixed topological type in $M$ with a free boundary on $\partial M$, assuming that $\partial M$ is strictly convex with respect to the inward unit normal. In this paper, we show that the strict convexity condition on $\partial M$ cannot be relaxed.
\end{abstract}

\maketitle

\section{Introduction}

Choi and Schoen \cite{Choi&Schoen} proved that for a $3$-dimensional Riemannian manifold $M$ with positive Ricci curvature, the space of embedded minimal surfaces in $M$ of a fixed genus is compact.
Hsieh and Wang \cite{Hsieh&Wang} have shown that the condition
of positive Ricci curvature cannot be relaxed.
In particular, they showed in \cite{Hsieh&Wang} that
the embedded minimal tori in $S^2\times S^1$ do not form a compact family,
for $S^2\times S^1$ equipped with the standard product metric
has nonnegative Ricci curvature.

On the other hand, Fraser and Li  \cite{Fraser&Li} proved the corresponding result of Choi and Schoen
for the free boundary minimal surfaces. In particular, they proved the following:

\begin{theorem}[Fraser-Li]\label{thm1}
Let $M$ be a compact $3$-dimensional Riemannian manifold with nonempty boundary $\partial M$.
Suppose $M$ has nonnegative Ricci curvature and $\partial M$ is strictly convex with respect to the inward unit normal.
Then the space of compact properly embedded minimal surfaces of fixed topological type in $M$ with free boundary in
$\partial M$ is compact in the $C^k$ topology for any $k\geq 2$.
\end{theorem}

Note that the existence of such a free boundary minimal annulus in $M$ was obtained by Maximo, Nunes and Smith \cite{Maximo&Nunes&Smith}. Inspired by the result of Hsieh and Wang \cite{Hsieh&Wang} mentioned above,
we show that the condition that $\partial M$ is strictly convex
could not be relaxed in Theorem \ref{thm1}.
To state our main theorem,
let $S^2_+$ be the $2$-dimensional hemisphere.
Equipped with the standard product metric, the product manifold $S^2_+\times S^1$  has nonnegative Ricci curvature,
and its boundary $\partial(S^2_+\times S^1)$
is totally geodesic and is not strictly convex.
The following is our main theorem.

\begin{theorem}\label{thm2}
Let $S^2_+\times S^1$ be equipped with the standard product metric.
Then there exists a sequence of compact {properly embedded} minimal annuli in $S^2_+\times S^1$ with free boundary in $\partial(S^2_+\times S^1)$ which
is not compact.
\end{theorem}

\section{Proof of the main theorem}

In this section, we prove Theorem \ref{thm2}.
Let $S^2_+$ be the $2$-dimensional hemisphere given by
$$S^2_+=\{(x,y,z)\in\mathbb{R}^3: x^2+y^2+z^2=1, y\geq 0\},$$
which is a $2$-dimensional manifold with boundary
$$\partial S^2_+=\{(x,0,z)\in\mathbb{R}^3: x^2+z^2=1\}.$$
Then $S^2_+$ has a local parametrization given by
\begin{equation}\label{0.9}
\overrightarrow{x}(r,\theta)=(\sin r\cos \theta,\sin r\sin\theta, \cos r),
\end{equation}
where $0\leq r\leq \pi$ and $0\leq \theta\leq \pi$.
Note that the boundary $\partial S^2_+$ corresponds to
the points when $\theta=0$ or $\theta=\pi$.
The induced metric
on $S^2_+$ by $\mathbb{R}^3$ is given by
$$dr^2+\sin^2 rd\theta^2.$$
Then the product space $S^2_+\times S^1$ is a $3$-dimensional manifold with boundary
$\partial (S^2_+\times S^1)=\partial S^2_+\times S^1$. Moreover,
\begin{equation}\label{0.1}
g=dr^2+\sin^2 rd\theta^2+dt^2
\end{equation}
is a Riemannian
metric on $S^2_+\times S^1$.
Equipped with the metric $g$ in (\ref{0.1}),
the boundary $\partial (S^2_+\times S^1)$ is a totally geodesic and hence is not strictly convex. It follows from (\ref{0.1}) that
\begin{equation}\label{0.2}
\begin{split}
\nabla_{\frac{\partial}{\partial\theta}}\frac{\partial}{\partial r}=\cot r\frac{\partial}{\partial\theta},
~~\nabla_{\frac{\partial}{\partial\theta}}\frac{\partial}{\partial t}=0,\\
\nabla_{\frac{\partial}{\partial r}}\frac{\partial}{\partial r}=0,~~
\nabla_{\frac{\partial}{\partial r}}\frac{\partial}{\partial t}=\nabla_{\frac{\partial}{\partial t}}\frac{\partial}{\partial r}=0,~~
\nabla_{\frac{\partial}{\partial t}}\frac{\partial}{\partial t}=0.
\end{split}
\end{equation}

Suppose $r(t)$ is a smooth function on $S^1$ such that
\begin{equation}\label{0.5}
r(t+2\pi)=r(t).
\end{equation}
Define
$f:[0,\pi]\times S^1\to S^2_+\times S^1$ given by
\begin{equation}\label{0.8}
f(\theta,t)=(\overrightarrow{x}(r(t),\theta),t),
\end{equation}
where $\overrightarrow{x}$ is defined as in (\ref{0.9}).
Then
$$\vec{e}_1:=\frac{1}{\sin r}\frac{\partial}{\partial\theta} ~\text{and}~~  \vec{e}_2:=\frac{1}{\sqrt{1+r'^2}}\left(r'\frac{\partial}{\partial r}+\frac{\partial}{\partial t}\right)$$
are unit vectors with respect to $g$,
and $\{\vec{e}_1,\vec{e}_2\}$ spans the tangent space of $[0,\pi]\times S^1$.
Hence the unit normal $\vec{n}$ is given by
$$\vec{n}=\frac{1}{\sqrt{1+r'^2}}\left(\frac{\partial}{\partial r}-r'\frac{\partial}{\partial t}\right).$$
By (\ref{0.2}), we see that
\begin{equation}\label{0.3}
(\nabla_{\vec{e}_1}\vec{n})^T=\nabla_{\vec{e}_1}\vec{n}=\left(\frac{\cot r}{\sqrt{1+r'^2}}\right) \vec{e}_1~~\mbox{ and }~~
(\nabla_{\vec{e}_2}\vec{n})^T=-\frac{r''}{(1+r'^2)^{\frac{3}{2}}} \vec{e}_2,
\end{equation}
where $(\nabla_{\vec{e}_2}\vec{n})^T$ denotes the tangential component of $\nabla_{\vec{e}_2}\vec{n}$.
It follows from (\ref{0.3}) that the principal curvatures
of $[0,\pi]\times S^1$ are $\frac{\cot r}{\sqrt{1+r'^2}}$ and $-\frac{r''}{(1+r'^2)^{\frac{3}{2}}}$,
which implies that the mean curvature is
$$H=\frac{1}{2\sqrt{1+r'^2}}\left(\cot r-\frac{r''}{1+r'^2}\right).$$
Hence, $[0,\pi]\times S^1$ is minimal if
\begin{equation}\label{0.4}
\frac{r''}{1+r'^2}=\cot r.
\end{equation}
On the other hand, the outward conormal vector $\vec{\nu}$
is the unique unit vector in $T_p([0,\pi]\times S^1)$ which is orthogonal to
$T_p\big(\partial ([0,\pi]\times S^1)\big)$. Thus, we have
\begin{equation}\label{0.6}
\vec{\nu}=\pm\frac{1}{\sin r}\frac{\partial}{\partial\theta}.
\end{equation}
Note that the tangent plane of $\partial(S^2_+\times S^1)$
is spanned by $\displaystyle\left\{\frac{\partial}{\partial r}, \frac{\partial}{\partial t}\right\}$.
Hence, it follows from (\ref{0.6}) that $\vec{\nu}$ is orthogonal to
$\partial(S^2_+\times S^1)$, i.e.,
\begin{equation}\label{0.7}
\vec{\nu}\perp \partial(S^2_+\times S^1).
\end{equation}

In view of (\ref{0.7}), we conclude that if $r(t)$ is a solution to (\ref{0.5}) and (\ref{0.4}),
then the embedding $f$ of $[0,\pi]\times S^1$ defined in (\ref{0.8})
is a free boundary minimal surface in $S^2_+\times S^1$.
As pointed out in \cite{Hsieh&Wang},
the general solutions
to (\ref{0.5}) and (\ref{0.4})
are elliptic functions for a sequence of initial values
$r(t=0)=c_1>c_2>c_3>\cdots\to 0$.
The corresponding sequence of embedded free boundary minimal surfaces certainly cannot have a convergent subsequence. This proves our main theorem.

%\section*{Acknowlegement}
%The first author was supported by the National Science and Technology Council (NSTC), Taiwan, with grant Number: 112-2115-M-032 -006 -MY2. The second author was supported by the National Research Foundation of Korea (NRF-2020R1A2C1A01005698 and NRF-2021R1A4A1032418). The third author was supported by the National Research Foundation of Korea (NRF-2021R1A2C1003365).

\bibliographystyle{amsplain}

\end{document}